\documentclass[twoside,10pt,fleqn]{article}
\usepackage{amsmath}
\usepackage{amsthm}
\usepackage{indentfirst,amsmath,amsfonts,amssymb,amsthm,cite}
\usepackage{epstopdf}
\usepackage{longtable,booktabs}
\usepackage{fancyhdr}   
\usepackage[dvips]{graphics}
\usepackage{graphicx}
\numberwithin{equation}{section}
\usepackage{color}
\usepackage{CJK}
\usepackage{float}
\usepackage{subfigure}

\begin{document}
\setcounter{page}{1}
\vspace*{0.2cm}  
\begin{center}
\baselineskip=7mm
{\LARGE\bf {A Markovian regime-switching stochastic SEQIR epidemic model with governmental policy}$^{\dagger}$}\\[1cm]
\normalsize\bf
Hongjie Fan$^{1,2,\ddagger}$, Kai Wang$^{2,\ddagger,*}$, Yanling Zhu$^{2,*}$ \\[5mm]
 {\footnotesize\it
 1. School of Economics and Management, Shandong Agriculture and Engineering University,\\
 Jinan 250100, Shandong, China\\[3mm]
 2.
 School of Statistics and Applied
Mathematics, Anhui University of Finance and Economics,\\ Bengbu
233030, Anhui, China}
\begin{figure}[b]
\footnotesize
\rule[-2.truemm]{5cm}{0.1truemm}\\[1mm]\baselineskip=5.mm
{$^\dagger$ Supported by the Natural Science Foundation of Anhui Province(no.2108085MA04), and the Social Science Foundation of the Ministry of Education of China(no.21YJAZH081).\\
$^{\ddagger}$ These authors have contributed equally to this paper as the first authors.\\
$^*$ Corresponding authors.\\
\vspace{0.05cm}
Emails: fanhj1997@163.com(HJ Fan); wangkai050318@163.com(K Wang); zhuyanling99@126.com(YL Zhu).}
\end{figure}
\end{center}
\noindent\hrulefill \newline
\begin{center}
\begin{minipage}{13cm}\baselineskip=7.mm
\vspace{-3mm} \noindent{\bf Abstract:}\ \
In this paper, a  stochastic SEQIR epidemic model with Markovian regime-switching is proposed and investigated.
The governmental policy and implement efficiency  are concerned by a generalized incidence function of the susceptible class.
We have the existence and uniqueness of the globally  positive solution to the stochastic model by using the Lyapunov method.
In addition,  we study the dynamical behaviors of the disease, and the sufficient conditions for the extinction  and persistence in mean  are obtained.
Finally,   numerical simulations are introduced to demonstrate the theoretical results.
\\[0.5cm]
\\[0.5cm]
\end{minipage}
\end{center}
\noindent \hrulefill \newline
\newpage
\section{Introduction }
In recent years,  COVID-19  has brought great disasters to people around the world.
It caused a lot of people lost their lives, but also lost a lot of natural resources.
For virus-borne diseases, in fact, people have been doing research for a long time.
The mathematical model also plays an indelible role in analyzing the spread dynamics of infectious diseases, and the most famous is SIR model, which was firstly studied by  Kermack and McKendrick in 1927 \cite{a9}.

After that, a  lot of models have been studied, such as SIS, SEIR, SIQR etc \cite{a5,a6,a7,a8,a10},
and researchers started to investigate the varying population \cite{c1,c2}.
Due to the heterogeneity and variation of the virus, the phenomenon of recurrence and reinfection appeared, and corresponding theoretical models were established \cite{b4,b5}.
And many factors are considered in the theoretical models, such as the vaccination, the imperfect quarantine and  media effects.
The purpose of authors is introduce the model that incorporates a partially effective quarantine policy to test the mathematical robustness of prior mathematical results to variations in  quarantine effectiveness \cite{b1,b2,b3}.
The media effect is also a important factor that influence the dynamical behaviors of the diseases \cite{b4,b5,b6}.
Due to the spread of information such as media, it will increase people's awareness of the disease and thus reduce infection.
Moreover, the governmental    prevention and control policy plays an important role in preventing the transmission of the epidemic.
 For example, Mandal et al \cite{ref2} proposed the model
\begin{equation}\label{wk11}
\aligned
\frac{dS}{dt}&= A-\beta(1-\rho_1)(1-\rho_2)SE   +b_1Q  -\xi S-pSM,\\
\frac{dE}{dt}&=\beta(1-\rho_1)(1-\rho_2)SE -(b_2+\alpha+\sigma+\xi)E,\\
\frac{dQ}{dt}&=b_2E-(b_1+c+\xi)Q ,\\
\frac{dI}{dt}&= \alpha E+cQ-(\eta+\xi+\delta)I  ,\\
\frac{dR}{dt}&= \eta I+\sigma E-\xi R+pSM ,
\endaligned
\end{equation}
with  $S(0), E(0), Q(0), I(0), R(0)$
are nonnegative constants.
 $S$ is susceptible individuals, $E$ is exposed individuals, $Q$ is quarantined individuals, $I$ is infected individuals, $R$ is recovered individuals.
Other parameters are described as Table \ref{tab: Parameters}. 
\begin{center}
\begin{longtable}{c c |c  c}
\caption{Description of parameters to model (\ref{wk11}).}
\label{tab: Parameters} \\
\toprule
Parameter  & Description & Parameter  & Description \\
\midrule
$A$ & the recruitment rate of $S$  &$b_1$  & the rate of $Q$ to $S$ \\
$\beta$ & the disease transmission  rate &  $c$   &the rate of $Q$  to $I$ \\
$\rho_1 \in (0,1)$ & the rate of $S$ with proper precaution  & $\xi$ & the natural death rate \\
$\rho_2 \in (0,1) $&  the rate of $E$ with proper precaution   &$\delta$  & the death rate due to  COVID-19  \\
$\alpha$  & the rate of $E$ to $I$  & $M$ &  governmental policies   \\
 $b_2$  &the rate of $E$ to $Q$  &$p$ &the rate of policies  implemented\\
 $\sigma$ &    the recover rate of $E$            &$\eta $      &   the recover rate of $I$\\
\bottomrule
\end{longtable}
\end{center}
We  know that there are many uncertainties   that may affect the dynamical behaviors of infectious diseases, such as climate, population movement, etc \cite{a1,a11}.
Therefore, many researchers studied the epidemic model with the white noise \cite{a5,a6,a8,a10}.
Authors of  \cite{a10} studied modified model (\ref{wk11}) by introducing the fluctuation   in the parameter $\beta$, so that $\beta\rightarrow\beta +\sigma_0 \frac{dB(t)}{dt}$.
In real eco-systems, population dynamics are often affected by random switching of the external environment.
For example, the disease transmission rate $\beta$ in the epidemic model is corrected for meteorological factors, because many viruses and bacteria have better survival and infectivity in humid, less ultraviolet conditions.
The SIS epidemic model with Markovian switching was firstly studied   in 2012 \cite{c3}.
Then  
the deterministic and stochastic SEQIR model was studied \cite{ref2,c4}. 
Based on the above work   and model (\ref{wk11}), we propose the following stochastic SEQIR model with    government policy  
\begin{equation}\label{wk12}
\aligned
\frac{dS}{dt}=&A(r(t))-\beta(r(t))[1-\rho_1(r(t))][1-\rho_2(r(t))]SE   +b_1(r(t))Q -\xi(r(t)) S-p(r(t))M(r(t)) h(S) \\
&-\sigma_0(r(t)) [1-\rho_1(r(t))][1-\rho_2(r(t))]SE \frac{dB(t)}{dt},\\
\frac{dE}{dt}=&\beta(r(t))[1-\rho_1(r(t))][1-\rho_2(r(t))]SE -[b_2(r(t))+\alpha(r(t))+\sigma(r(t))+\xi(r(t))]E dt\\
&+\sigma_0(r(t)) [1-\rho_1(r(t))][1-\rho_2(r(t))]SE \frac{dB(t)}{dt},\\
\frac{dQ}{dt}=&b_2(r(t))E-[b_1(r(t))+c(r(t))+\xi(r(t))]Q ,\\
\frac{dI}{dt}=&\alpha(r(t)) E+c(r(t))Q-[\eta(r(t))+\xi(r(t))+\delta(r(t))]I  ,\\
\frac{dR}{dt}=&\eta(r(t)) I+\sigma(r(t)) E-\xi (r(t)) R+p(r(t)) M(r(t))h(S),
\endaligned
\end{equation}
where $B(t)$ is the standard Brownian motion  with the density of noises $\sigma_0^2$.
And   $pM h(S)$ is  a twice-differentiable  function of governmental policies, which enables the model to be more realistic.
That implies, not only can government policies change, but they can vary depending on the number of susceptible people.
We assume that $h(0)=0$ and   $0\leq h(S) \leq S h'(0)$ for all $S\geq 0$, see reference \cite{c1} for  more details.

The remainder of this study is structured as follows:
In Section 2, we show some preliminaries that will be used in the following sections.
In Section 3, we devote to Existence and uniqueness of the positive solution to the stochastic model (\ref{wk12}).
In Section 4, we obtain the condition of extinction of the disease.
In Section 5, the condition for the persistence in the mean  of the disease is obtained.
In Section 6, we carry out the numerical simulation of  and stochastic model.
This paper ends with  conclusion and discussion in Section 7.

\section{Preliminaries}

Let $(\Omega, \mathfrak{F}, \{\mathfrak{F}_t \}_{t\geq0}, \mathbb{P})$ be a complete probability space with a filtration $ \{\mathfrak{F}_t \}_{t\geq0}$ satisfying the usual conditions
(i.e., it is right continuous and $\mathfrak{F}_0 $ contains all $\mathbb{P}$-null sets), $B(t)$ be defined on this probability space, and $\mathbb{R}_+^n = \{(x_1,...x_n)\in \mathbb{R}^n : x_i \geq 0 ,i = 1,2,...,n \}$.
Let $r(t),\ t\geq 0$ be a right continuous  Markov chain on the complete probability space $(\Omega, \mathfrak{F}, \{\mathfrak{F}_t \}_{t\geq0}, \mathbb{P})$, taking values in a finite state space $\mathbb{S}={1,...,N}$ and with infinitesimal generator $\Gamma=(q_{ij})\in \mathbb{R}^{N\times N}$.
That is $r(t)$ satisfies
\[\aligned
P(r(t+\Delta)=j | r(t)=i)=\left\{
         \begin{array}{ll}
           q_{ij}\Delta+o(\Delta), & if\  i \neq j, \\
           1 +q_{ij}\Delta+o(\Delta), & if\  i = j,
         \end{array}
       \right.
       \endaligned\]
where $\Delta >0$ and $q_{ij}>0$ is the transition rate from $i$ to $j$ if $i\neq j$ while $q_{ii}=-\sum_{i\neq j} q_{ij}$ for ech $i\in \mathbb{S}$.
That implies $r_t$ is irreducible and has  a unique stationary distribution $\pi =(\pi_1,...,\pi_N)$, which can be determined by $\pi \Gamma=0$, subjected to $\sum_{k=1}^{N}\pi_k=1 $, $\pi_k>0$, for any $k\in \mathbb{S}$. And assume that the right continuous Markov chain $r(t)$ is independent of the Brownian motion $B(t)$ for $t\geq 0$.

In general, considering the $n$-dimensional stochastic differential equation
\begin{equation*}\label{wk21}
\aligned
dx(t)=F(x(t),r(t))dt+g(x(t),r(t))dB(t),\ \text{for}\  t\geq 0, \\
\endaligned
\end{equation*}
where $F(x(t),r(t))$ is  defined on $\mathbb{R}^{n} \times  \mathbb{S} \rightarrow \mathbb{R}^{n}  $
and $g(x(t),r(t))$ is an $n \times m$ matrix, which are locally Lipschitz functions in $x$.
$B(t)$  and $r(t)$ are $m$-dimensional Brownian motion and the right continuous Markov chain in the above discussion.

For each $k\in \mathbb{S}$,
  a given function
  $
V (x,k): \mathbb{R}^{n} \times \mathbb{S} \rightarrow \mathbb{R}^{n} $,
such that $V(x,k)$ is twice continuously differential with respect to the first variable $x$,
we denote
\[\aligned
LV(x,i) = V_{t}(x,i)+V_{x}(x,i) F(x,i) + \frac{1}{2}trace\left[ g^{T}(x,i)V_{xx}(x,i) g(x,i) \right]  +\sum_{k\in \mathbb{S}} q_{ik}V(x,k),
\endaligned
\]
where $
V_{t}(x,k)=\frac{\partial V(x,k)}{\partial t},
V_{x}(x,k)=(\frac{\partial V(x,k)}{\partial x_{1}}, \ldots, \frac{\partial V(x,k)}{\partial x_{n}}),
V_{x x}(x,k)
=(\frac{\partial^{2} V(x,k)}{\partial x_{i} \partial x_{j}})_{n \times n}$.
By It$\hat{o}$'s formula we have
\begin{equation*}\label{wk22}
\aligned
dV(x,r(t)) = LV(x,r(t))dt +V_{x}(x,r(t))g(x,r(t))dB(t).
\endaligned
\end{equation*}

For simplicity, we define $\hat{\varrho}=\min_{k\in \mathbb{S}} {\varrho(k)}$ and $\check{\varrho}=\max_{k\in \mathbb{S}} {\varrho(k)}$.

\section{Existence and uniqueness of the positive solution to the stochastic model}
In this section, we discuss the existence and uniqueness of the positive solution   in the stochastic model (\ref{wk12}).
First of all, we give a useful lemma.
\\[3mm]
\textbf{Lemma 3.1.}
Let 
  $\widetilde{N}(t)=S(t)+E(t)+ Q(t)+I(t)+R(t)$, then we can obtain the positively invariant set
\[\aligned\Omega=\left\{ ( S(t),E(t), Q(t),I(t),R(t))\in \mathbb{R}_{+}^{5}:\   \frac{\hat{A}}{\check{\xi} +\check{\delta} } \leq S(t)+E(t)+ Q(t)+I(t)+R(t) \leq \frac{\check{A}}{\hat{\xi} }\right\}.\endaligned\]
\\
\textbf{Proof.}
The proof of the Lemma is similar with   Remark 3.1 of \cite{a1}, therefore it is omitted.
$\square$\\
\\[3mm]
\textbf{Theorem 3.1.} For $r_0 \in \mathbb{S} $ and any initial value $(S(r_0), E(r_0), Q(r_0), I(r_0), R(r_0)) \in \mathbb{R}_{+}^{5}$,
there is a unique positive solution $(S, E, Q, I, R)$
of model {(\ref{wk12})} for  $t \geq 0$, which will remain in $\mathbb{R}_{+}^{5}$ with probability one.
\\[3mm]
\textbf{Proof.}
Since the coefficients of model  (\ref{wk12}) are locally Lipschitz continuous on $\mathbb{R}_{+}$,
then for any $r_0 \in \mathbb{S} $ and initial value $(S(r_0), E(r_0), Q(r_0), I(r_0), R(r_0)) \in \mathbb{R}_{+}^{5}$,
 there is a unique local solution $(S(t), E(t), Q(t), I(t), R(t))$ to model {(\ref{wk12})} on $[0,\tau_{e}]$.
 We need to prove  $\tau_{e}=+\infty$ a.s. to show that the solution exists globally.
 Let $k_{0}\geq1$ be sufficiently large such
that $S, E, Q, I, R$  are lie in $[k_{0},\frac{1}{k_{0}}]$.
Let $k\geq k_{0}$, define the stopping time
\[\aligned  \tau_{k}=\inf \left\{t \in\left[0, \tau_{e}\right): S(t) \notin\left(\frac{1}{k}, k\right) \text { or } E(t) \notin\left(\frac{1}{k}, k\right)
\text { or } Q(t) \notin\left(\frac{1}{k}, k\right) \text { or } I(t) \notin\left(\frac{1}{k}, k\right) \text { or } R(t) \notin\left(\frac{1}{k}, k\right)\right\}.
\endaligned\]
Let $\tau_{+\infty}=\lim _{k \rightarrow+\infty} \tau_{k}$, whence $\tau_{+\infty} \leq \tau_{e}$. If we can show $\tau_{+\infty}=+\infty$,
then $\tau_{e}=+\infty$ and  model (\ref{wk12}) has a unique solution $(S, E, Q, I, R)$ for  $t\geq 0$.

Obviously, we only need to show
$\tau_{+\infty}=+\infty$.
If the assertion is false, then there exist constants $T>0$ and $\epsilon \in (0,1)$ such that
$\mathbb{P}(\tau_{+\infty}
\leq T) \geq \epsilon,$ which yields that
there is an integer $k_{1} \geq k_{0}$ such that
$ \mathbb{P}(\tau_{k} \leq T) \geq \epsilon,\ \text{for\ }k \geq k_{1}. $
Define a Lyapunov function $V:\mathbb{R}_+^5 \rightarrow \mathbb{R}_+$:
 \[\aligned
V(S,E,Q,I,R)=S-1-\ln S+E-1-\ln E+Q-1-\ln Q+I-1-\ln I+R-1-\ln R.
\endaligned \]
By It\^{o}'s formula
\[\aligned
LV
\leq & A(r(t))+ \beta(r(t))(1-\rho_1(r(t)))(1-\rho_2(r(t)))E+5\xi(r(t))+b_1(r(t))+b_2(r(t))+c(r(t))+\alpha(r(t))\\
&+\sigma(r(t))+\eta(r(t))
+\delta(r(t))+ \frac{1}{2}\sigma_0^2(r(t)) (1-\rho_1(r(t)))^2(1-\rho_2(r(t)))^2(E^2+S^2)+p(r(t))M(r(t)) \frac{h(S)}{S} \\
\leq & \hat{A}+\frac{\hat{A}\hat{\beta} (1-\check{\rho}_1) (1-\check{\rho}_2)}{\check{\xi}}+4\hat{\xi}+\hat{b_1}+\hat{b_2}+\hat{\alpha}+ \hat{\sigma}+ \hat{\xi} +\hat{c}+\hat{\eta}+\hat{\delta}+\hat{p}\hat{M}h'(0)+\hat{\sigma}_0^2 (1-\check{\rho}_1)^2(1-\check{\rho}_2)^2
   \frac{\hat{A}}{\check{\xi} }^2.
\endaligned \]
The remainder of the proof follows that in \cite{a2}, and hence it is omitted here.
$\square$

\section{Extinction }
In this section, we obtain the condition for the extinction of the disease.
We define
 \[\aligned  R_s^*=
\frac{\sum_{k=1}^{N}\pi(k) \beta(k) w_1(k) \frac{\check{A}}{\hat{\xi}} }
{\sum_{k=1}^{N}  \pi(k)  \left[w_2(k)+ \frac{\sigma_0^2(k)}{2 } w_1^2(k)  \left(\frac{\check{A}}{\hat{\xi}} \right)^2\right]},
\endaligned \]
where $w_1(k)=(1-\rho_1(k))(1-\rho_2(k))$, and $w_2(k)=b_2(k)+\alpha(k)+\sigma(k)+\xi(k)$ for $k\in \mathbb{S}$.
\\[3mm]
\textbf{Theorem 4.1.} 
If  $  \beta(k) \geq \sigma_0^2(k)w_1(k)  \frac{\check{A} }{\hat{\xi} }$ and $R_s^* <1$ for any $k\in \mathbb{S}$, then, for any $r_0 \in \mathbb{S} $ and the initial value $(S(r_0), E(r_0), Q(r_0), I(r_0), R(r_0)) \in \mathbb{R}_{+}^{5}$, the disease will go to extinction exponentially with probability one, i.e.,
 \[\aligned \limsup_{t\rightarrow +\infty}E(t)=\limsup_{t\rightarrow +\infty}Q(t)=\limsup_{t\rightarrow +\infty}I(t)=0\ a.s..
\endaligned \]
\\
\textbf{Proof.} 
By applying It\^{o}'s formula to the second equation of  model (\ref{wk12}), we have
\begin{equation}\label{30}
\aligned
 d\ln E(t)=& \left[\beta(r(t)) (1-\rho_1(r(t)))(1-\rho_2(r(t))) S(t) -(b_2(r(t))+\alpha(r(t))+\sigma(r(t))+\xi(r(t)))\right. \\
 & \left. -\frac{1}{2}\sigma_0^2(r(t)) (1-\rho_1(r(t)))^2(1-\rho_2(r(t)))^2 S^2(t) \right] dt\\
 &+\sigma_0(r(t))(1-\rho_1(k))(1-\rho_2(k))S(t)dB(t).
\endaligned
\end{equation}
Integrating both sides of the above equation form $0$ to $t$ gives
 \[\aligned
 \ln E(t) -\ln e_0
 =&
 \int_0^t \left[\beta(r(u)) w_1(r(u))S(u) -\frac{1}{2}\sigma_0^2(r(u)) w_1^2(r(u)) S^2(u) -w_2(r(u)) \right] du
 +\Upsilon(t),
\endaligned \]
where
$$
\Upsilon(t)=\int_0^t \sigma_0(r(u)) w_1(r(u)) S(u)dB(u).
$$
 By applying  of the large number theorem for local continuous  martingales, we can get
 \[\aligned
  \frac{\ln E(t) -\ln e_0}{t}=\frac{1}{t} \int_0^t \left[\beta(r(u)) w_1(r(u))S(u) -\frac{1}{2}\sigma_0^2(r(u)) w_1^2(r(u)) S^2(u) -w_2(k) \right] du.
  \endaligned \]
Further, for any $k\in $ $\mathbb{S}$, we define the function $x\longmapsto F(x)= \beta(k) w_1(k) x -\frac{1}{2}\sigma_0^2(k) w_1^2(k) x^2 -w_2(k) $.
By the positively invariant set $\Omega$ and the assumption $  \beta(k) \geq \sigma_0^2(k)w_1(k)  \frac{\check{A} }{\hat{\xi} }$,
we know $F(x)$ is increasing on $\left[0,\beta(k)/ \sigma_0^2(k) w_1(k)\right]$, 
then
  \[\aligned
  \frac{\ln E(t) -\ln e_0}{t}  \leq
  \frac{1}{t} \int_0^t \left[\beta(r(u)) w_1(r(u))\frac{\check{A}}{\hat{\xi}} -\frac{1}{2}\sigma_0^2(r(u)) w_1^2(r(u)) \frac{\check{A}^2}{\hat{\xi}^2} -w_2(r(u)) \right] du.
 \endaligned \]
Taking the superior limit on both sides of the above equation, we can obtain
 \[\aligned
 \limsup_{t \rightarrow \infty} \frac{\ln E(t) -\ln e_0}{t}
 &\leq  \limsup_{t \rightarrow \infty}
 \frac{1}{t} \int_0^t \left[\beta(r(u)) w_1(r(u))\frac{\check{A}}{\hat{\xi}} -\frac{1}{2}\sigma_0^2(r(u)) w_1^2(r(u)) \frac{\check{A}^2}{\hat{\xi}^2} -w_2(r(u)) \right] du
  \endaligned \]
   It follows from the ergodic property of the Markov chain that
  \[\aligned
 \limsup_{t \rightarrow \infty} \frac{\ln E(t) -\ln e_0}{t}
 & \leq  \sum_{k=1}^{N} \pi (k)
\left[\beta(k) w_1(k)\frac{\check{A}}{\hat{\xi}} -\frac{1}{2}\sigma_0^2(k) w_1^2(k) \frac{\check{A}^2}{\hat{\xi}^2} -w_2(k) \right]\\
&\leq  \sum_{k=1}^{N} \pi (k) \left (w_2(k) + \frac{1}{2}\sigma_0^2(k) w_1^2(k) \frac{\check{A}^2}{\hat{\xi}^2}  \right )
 \left(R_s^*-1  \right)   <0 \ a.s.,
 \endaligned \]
which implies that
$  \limsup_{t\rightarrow +\infty} E(t)=0 \ a.s.. 
$
Thus for any constant $\epsilon>0$, there exists a positive constant $T$ such that $E(t,\omega)\leq \epsilon$ for $ t>T$, which together with the third equation of model (\ref{wk12}) yields
 \[\aligned
\frac{dQ(t)}{dt}\leq \check{b}_2E(t)-(\hat{b}_1+\hat{c}+\hat{\xi})Q(t)\leq \check{b}_2\epsilon -(\hat{b}_1+\hat{c}+\hat{\xi})Q(t).
\endaligned\]
By using the comparative theorem, we have
 \[\aligned
\limsup_{t\rightarrow +\infty}Q(t) \leq \frac{\check{b}_2\epsilon}{\hat{b}_1+\hat{c}+\hat{\xi}} \text{\ a.s..}
\endaligned\]
Letting $\epsilon\rightarrow 0$, $\limsup_{t\rightarrow +\infty}Q(t)=0$ a.s..
Similarly,  we can obtain that $\limsup_{t\rightarrow +\infty}I(t)=0$ a.s..
$\square$

\section{Persistence in mean}
In this section, we obtain the condition for the persistence of the disease.
We define
\[\aligned \widetilde{R}_s^*= \frac{     \sum_{k=1}^{N }\pi(k)
\beta(k)w_1(k) \frac{\check{A}}{\hat{\xi}}
 }{ \sum_{k=1}^{N }\pi(k) \left [ w_2(k)+\psi_1(k) + \frac{\sigma_0^2(k)}{2}w_1^2(k) \left(\frac{\check{A}}{\hat{\xi}} \right)^2  \right]
 },
 \endaligned\]
 where $\psi_1 $ is  determined by equation (\ref{411}).
 \\[3mm]
\textbf{Theorem 5.1.} 
If $  \beta(k) \geq \sigma_0^2(k)w_1(k)  \frac{\check{A} }{\hat{\xi} }$ and $\widetilde{R}_s^* >1$ for any $k\in\mathbb{S}$,
then, for any $r_0 \in \mathbb{S} $ and the initial value $(S(r_0), E(r_0), Q(r_0), I(r_0), R(r_0)) \in \mathbb{R}_{+}^{5}$, the disease will persist in mean, and the solution $(S(t), E(t), Q(t), I(t), R(t)) $ to  model (\ref{wk12}) has the following properties
\begin{equation}\label{46}\aligned
\liminf_{t\rightarrow +\infty}  \frac{1}{t}\int_{0}^{t}  E(s) ds  \geq
\frac{\Lambda}{\sum_{k=1}^{N} \pi(k) \psi_2(k)} (\widetilde{R}_s^*-1)\ a.s..
\endaligned\end{equation}
\begin{equation}\label{47}
\aligned
\liminf_{t\rightarrow +\infty}  \frac{1}{t}\int_{0}^{t}  Q(s) ds  \geq
\frac{\hat{b}_2  \Lambda}{(\check{b}_1+\check{c}+\check{\xi}) \sum_{k=1}^{N} \pi(k) \psi_2(k) }   (\widetilde{R}_s^*-1)\ a.s..
\endaligned\end{equation}
\begin{equation}\label{48}
\aligned
\liminf_{t\rightarrow +\infty}  \frac{1}{t}\int_{0}^{t}  I(s) ds  \geq
\frac{ 1}{\check{\eta} +\check{\xi} +\check{\delta}}  \left( \hat{\alpha} +\frac{\hat{c}\hat{b}_2  }{(\check{b}_1+\check{c}+\check{\xi})  }\right)
  \frac{\Lambda}{\sum_{k=1}^{N} \pi(k) \psi_2(k)}(\widetilde{R}_s^*-1)\ a.s.,
\endaligned\end{equation}
where $\psi_2 $ and $\Lambda$ are determined by equations (\ref{410}) and (\ref{49}), respectively.
 \\[3mm]
\textbf{Proof.} 
By   inequation (\ref{30}) of section 4, for any $k\in $ $\mathbb{S}$, we obtain that 
 \[\aligned
d \ln(E) = F\left(S,k\right)dt+\sigma_0(k)w_1(k)S(t)dB(t),
 \endaligned \]
where
\[\aligned
  F(S,k)= \beta(k) w_1(k) S -\frac{1}{2}\sigma_0^2(k)  w_1^2(k) S^2 - w_2(k).
\endaligned\]
Then we can obtain
\begin{equation}\label{40}
\aligned
F(S,k)
 & \geq  F\left(\frac{\check{A}}{\hat{\xi}},k \right)- \left( \beta(k) w_1(k)- \frac{\sigma_0^2(k)}{2} w_1^2(k) \frac{\check{A}}{\hat{\xi}}   \right)   \left(  \frac{\check{A}}{\hat{\xi}} -S \right) \\
& \geq  F\left(\frac{\check{A}}{\hat{\xi}},k \right)
-  \left( \check{\beta} w_1(k)-\frac{{\sigma}_0^2(k)}{2} w_1^2(k) \frac{\check{A}}{\hat{\xi}}  \right)  \left( \frac{\check{A} }{\hat{\xi}} -S\right).
\endaligned
\end{equation}
By the first equation of   model (\ref{wk12}), we have
\[\aligned
dS  \geq &
\left(  A(k)-(\beta(k) w_1(k) E+\xi(k))S - p(k)M(k)h(S)  \right) dt
-\sigma_0(k) w_1(k)SE dB(t)\\
\geq &
\left(  A(k)-(\check{\beta} w_1(k) E+\xi(k))S - p(k)M(k)h(S)  \right) dt
-\sigma_0(k) w_1(k)SE dB(t)\\
\geq &
\left[     {A(k)}\frac{  \hat{\xi}}{\check{A}}\left( \frac{\check{A}}{\hat{\xi}}-S \right)-{\xi(k)}S\left( 1-\frac{{A(k) } \hat{\xi}}{\check{A}{\xi(k)}}+\frac{{p(k)}{M(k)}h(S)}{\xi(k)S} \right)
-{\check{\beta}}w_1(k)\frac{\check{A}}{\hat{\xi}}E    \right] dt-\sigma_0(k) w_1(k)SE dB(t).
\endaligned\]
Then
\begin{equation}\label{41}
\aligned
-\left( \frac{\check{A}}{\hat{\xi}}-S \right)dt
\geq &  \frac{\check{A}}{{A(k)} \hat{\xi}} \left[
  - dS-{\xi(k)}S\left( 1-\frac{{A(k) } \hat{\xi}}{\check{A}{\xi(k)}}+\frac{{p(k)}{M(k)}h(S)}{\xi(k)S} \right)
-{\check{\beta}}w_1(k)\frac{\check{A}}{\hat{\xi}}E       \right] dt
-\frac{\check{A} \sigma_0(k) w_1(k)}{ {A(k)} \hat{\xi}}SE dB(t)\\
\geq &
\frac{\check{A}}{ {A(k)} \hat{\xi}} \left[
  - {\xi(k)} \frac{\check{A}}{\hat{\xi}}\left( 1-\frac{{A(k) } \hat{\xi}}{\check{A}{\xi(k)}}+\frac{{p(k)}{M(k)}h'(0)}{\xi(k)} \right)
  dt   \right] -\frac{\check{A}}{ {A(k)} \hat{\xi}} dS - \frac{\check{A}}{ {A(k)} \hat{\xi}}   \check{\beta}w_1(k)\frac{\check{A}}{\hat{\xi}}E dt                  \\
&-\frac{\check{A} \sigma_0(k) w_1(k)}{ {A(k)} \hat{\xi}}SE dB(t).
\endaligned
\end{equation}
Therefore, we can obtain
\begin{equation}\label{43}
\aligned
d \ln(E) = &
F\left(S,k \right)dt+\sigma_0(k)w_1(k)S(t)dB(t)\\
  \geq&  F\left(\frac{\check{A}}{\hat{\xi}},k \right)dt
-  \left( \check{\beta}w_1(k)-\frac{\hat{\sigma}_0^2}{2}w_1^2(k) \frac{\check{A}}{\hat{\xi}}  \right)  \left( \frac{\check{A} }{\hat{\xi}} -S\right)dt+\sigma_0(k)w_1(k)S(t)dB(t)\\
\geq&  F\left(\frac{\check{A}}{\hat{\xi}}, k \right)dt
- \psi_1(k) dt  -\psi_2(k) Edt- \psi_3(k) dS   + \varphi_1(k) SdB(t) + \varphi_2(k)SEdB(t),
\endaligned
\end{equation}
where
\begin{equation}\label{411}\aligned
\psi_1(k):=\left( \check{\beta}w_1(k)-\frac{\hat{\sigma}_0^2}{2}w_1^2(k) \frac{\check{A}}{\hat{\xi}}  \right)
\frac{\check{A}^2 {\xi(k)}}{ {A(k)} \hat{\xi}^2}
   \left( 1-\frac{{A(k) } \hat{\xi}}{\check{A}{\xi(k)}}+\frac{{p(k)}{M(k)}h'(0)}{\xi(k)} \right),
    \endaligned\end{equation}
\begin{equation}\label{410}\aligned
\psi_2(k):=\left( \check{\beta}w_1(k)-\frac{\hat{\sigma}_0^2}{2}w_1^2(k) \frac{\check{A}}{\hat{\xi}}  \right)
    \frac{\check{A}^2}{ {A(k)} \hat{\xi}^2}   \check{\beta}w_1(k),
    \endaligned\end{equation}
\[\aligned
\psi_3(k):= \left( \check{\beta}w_1(k)-\frac{\hat{\sigma}_0^2}{2}w_1^2(k) \frac{\check{A}}{\hat{\xi}}  \right)\frac{\check{A}}{ {A(k)} \hat{\xi}},
   \endaligned\]
\[\aligned
\varphi_1(k) :=& \sigma_0(k)w_1(k),
\endaligned\]
\[\aligned
\varphi_2(k) : = -\left( \check{\beta}w_1(k)-\frac{\hat{\sigma}_0^2}{2}w_1^2(k) \frac{\check{A}}{\hat{\xi}}  \right)\frac{\check{A} \sigma_0(k) w_1(k)}{ {A(k)} \hat{\xi}}.
   \endaligned\]
We define a Lyapunov function $U:\ \mathbb{R}_+ \times \mathbb{S} \rightarrow \mathbb{R} $, which is   $U(E,k)=\ln (E) + \omega(k)$. Then
\begin{equation}\label{44}
\aligned
dU\geq F\left(\frac{\check{A}}{\hat{\xi}}, k \right)dt
- \psi_1 (k)dt  -\psi_2 (k)Edt +\sum_{k=1}^{N}\gamma_{kl}\omega(l) dt- \psi_3(k) dS   +  \varphi_1 (k)SdB(t) + \varphi_2(k) SEdB(t)
   \endaligned
   \end{equation}
Since the generator matrix is  irreducible, for $P_0=(P(1),......P(N)) $ with $P(k)=F\left(\frac{\check{A}}{\hat{\xi}},k \right)$, there is a $\omega=(\omega(1),......\omega(N))$ satisfying the  Poisson system
\[\aligned
\Gamma \omega = \sum_{h=1}^{N} \pi_h P_0(h) \textbf{1} -P_0,
\endaligned\]
where \textbf{1} is  a unit of vector of $\mathbb{R}^N$.
That implies
\begin{equation}\label{wk45}
\aligned
F\left(\frac{\check{A}}{\hat{\xi}},k \right)+\sum_{k=1}^{N}\gamma_{kl}\omega(l)  =\sum_{k=1}^{N}\pi(k)F\left(\frac{\check{A}}{\hat{\xi}},k \right)
\endaligned
\end{equation}

Substituting the above equality (\ref{wk45}) into the inequality (\ref{44}), integrating from $0$ to $t$ and dividing by $t$ on both sides,
\[\aligned
\frac{U(t)-U(0)}{t}\geq &\sum_{k=1}^{N}F\left(\frac{\check{A}}{\hat{\xi}},k \right)
-\frac{1}{t}\int_{0}^{t} \psi_1 ds
-\frac{1}{t}\int_{0}^{t} \psi_2 E(s) ds
-\frac{1}{t}\int_{0}^{t} \psi_3 dS\\&
-\frac{1}{t}\int_{0}^{t} \varphi_1 S(s)  dB(s)
-\frac{1}{t}\int_{0}^{t} \varphi_2 S(s)E(s) dB(s)
\endaligned\]
Combined with the boundness of $S$ and the strong law of the large number theorem for continuous local martingales,
\[\aligned
\lim_{t\rightarrow +\infty }  \frac{1}{t}\int_{0}^{t} \psi_3(k) dS
+\lim_{t\rightarrow +\infty } \frac{1}{t}\int_{0}^{t} \varphi_1 (k)S(s)  dB(s)
+\lim_{t\rightarrow +\infty } \frac{1}{t}\int_{0}^{t} \varphi_2 (k)S(s)E(s) dB(s)=0\ a.s..
\endaligned\]
Moreover, we have $\lim_{t\rightarrow +\infty } \frac{U(t)-U(0)}{t} =0 \ a.s..  $
Therefore
\[\aligned
 \frac{1}{t}\int_{0}^{t} \psi_2(k) E(s) dt  \geq&
 \sum_{k=1}^{N}\pi(k) F\left(\frac{\check{A}}{\hat{\xi}},k \right)
-\frac{1}{t}\int_{0}^{t} \psi_1 (k)ds\\
\geq &
\sum_{k=1}^{N}\pi(k)  \left[
\beta(k)w_1(k) \frac{\check{A}}{\hat{\xi}} -\frac{\sigma_0^2(k)}{2}w_1^2(k) \left(\frac{\check{A}}{\hat{\xi}} \right)^2 -w_2(k)
\right]
-\sum_{k=1}^{N} \pi(k) \psi_1(k)\\
\geq &
\sum_{k=1}^{N}\pi(k)  \left[
   \frac{\sigma_0^2(k)}{2}w_1^2(k)\left(\frac{\check{A}}{\hat{\xi}} \right)^2 +w_2(k)+\psi_1(k)
\right](\widetilde{R}_s^*-1):=  \Lambda (\widetilde{R}_s^*-1),
\endaligned\]
where
\begin{equation}\label{49}\aligned
 \Lambda=\sum_{k=1}^{N}\pi(k)  \left[
   \frac{\sigma_0^2(k)}{2}w_1^2(k) \left(\frac{\check{A}}{\hat{\xi}} \right)^2 +w_2(k)+\psi_1(k)
\right].
\endaligned\end{equation}

Similarly, we can obtain the equation (\ref{46}), (\ref{47}) and  (\ref{48}). $\square$
\\[3mm]
\textbf{Remark 5.1.}
By the assumption of $h(S) $ and 
the condition of Theorem 5.1, we can   obtain
$1-\frac{{A(k) } \hat{\xi}}{\check{A}{\xi(k)}}+\frac{{p(k)}{M(k)}h'(0)}{\xi(k)}>0$ and 
$  \beta(k) \geq \frac{1}{2}\sigma_0^2(k) w_1(k) \frac{\check{A} }{\hat{\xi} }$ for $k\in \mathbb{S}$
then $\psi_i,\ i=1,2,3$ are nonnegative.
Therefore, we can have that $\widetilde{R}_s^*\leq  {R}_s^*$, and the equation holds if and only if $  \beta(k) = \frac{1}{2}\sigma_0^2(k) w_1(k) \frac{\check{A} }{\hat{\xi} }$ for any $k\in \mathbb{S}$, then ${R}_s^*$ is the basic reproduction number of   stochastic model (\ref{wk12}). Comparing ${R}_s^*$ with the basic reproduction number $R_0$ of   deterministic model (\ref{wk11}), they are same if model (\ref{wk12}) has no white noise and Markov chain $r(t)$ has only one state.
That implies   model (\ref{wk12}) we proposed is more generalized and it is  more suitable to the complex environment.

\section{Simulation}
In this section, we give two examples by Milstein’s Higher Order Method \cite{a3,a4}, and set Markov chain $r(t)$ by $\mathbb{S}=\{1,2,3,4\}$. Let  $\Delta=0.0001$ is the step and the generator $\Gamma$ is
\[\aligned
\Gamma=\left(
\begin{array}{cccc}
-10 & 3& 2 &5 \\
6 & - 9 & 2 &1\\
3 & 3 &-8 & 2\\
1 & 5 & 3 &-9
\end{array}
\right),
\endaligned\]
then
\[\aligned
P=e^{\Delta \Gamma}\left(
\begin{array}{cccc}
0.9990 &0.0003 &0.0002 &0.0005\\
 0.0006 &0.9991 &0.0002 &0.0001 \\
 0.0003 &0.0003 &0.9992 &0.0002 \\
 0.0001& 0.0005 &0.0003 &0.9991
\end{array}
\right).\endaligned\]
Similarly, we have the stationary distribution of $r(t)$, which is $\pi=(0.2622, 0.2879, 0.2227, 0.2272)$.
The simulation of $r(t)$ is shown as Figure 1.
\\[3mm]
\textbf{Example 1.}
 Let
$ A=(0.0008,0.0005,0.0070,0.0010)$,    
$\beta=(0.006,0.018,0.049,0.08)$,

$\xi=(0.011,0.010,0.019,0.02)$,
$b_1=(0.05,0.06,0.010,0.08)$,
$b_2=(0.05,0.04,0.06,0.07)$,

$c=(0.08,0.07,0.09,0.10)$,
$\sigma=(0.003,0.005,0.006,0.004)$,
$\rho_1=(0.001,0.005,0.010,0.009)$,

$\rho_2=(0.001,0.005,0.007,0.003)$,
$\alpha=(0.016,0.0015,0.0017,0.0019)$,
$p=(0.001,0.002,0.003,0.004)$,

$\eta=(0.02,0.018,0.019,0.0021)$,
$\delta=(0.05,0.06,0.04,0.08)$,
$M=(0.001,0.002,0.003,0.004)$,

$\sigma_0=(0.008,0.065,0.007,0.006) $.\\
In addition, the initial values of the system are $S(0)=20,E(0)=20, Q(0)=15,I(0)=10, R(0)=0$.

We can have $R_s^*=0.1277<1$, which follows theorem 4.1.
That implies the disease will go out eventually, as shown in Figure 2.
In Figure 3 it depicts $E(t),Q(t)$ and $I(t)$ under the states $r(t)=1,2,3,4$, respectively.
The disease of the subsystem  with the states $r(t)=2,3$ will go out,
which are different from the results of model (\ref{wk12}) with regime switching.
\\[3mm]
\textbf{Example 2.}
We only change that
$A=(0.70,0.245,0.890,0.41)$, 
$\beta=(0.016,0.018,0.019,0.008)$,  
then $\widetilde{R}_s^*=2.5861>1$ and simulations are shown as Figure 4.
It is follows from theorem 5.1 that  the disease will persist in mean.
However, that will be distinct in the subsystem with the state  $r(t)=2$  in Figure 5, which depict $E(t),Q(t)$ and $I(t)$ under $r(t)=1,2,3,4$, respectively.
\\[3mm]
\textbf{Remark 6.1.}
Then it can be see that the stochastic epidemic model with regime switching is influenced by several states,
but it's not determined by one state.
Therefore the model with regime switching is more realistic and suitable to describe the dynamics of the diseases in the complex and changing environment.

\section{Conclusion and discussion}
In this paper, we propose  a SEQIR epidemic model with  both white and telegraph noises,
and investigate the dynamic behaviors of the diseases.
The governmental policy  and the efficiency of policy implemented are both considered in the stochastic model by constructing the generalized function.
Firstly, we get the positively invariant set of the classes of the stochastic model,
and have the existence and uniqueness of the globally positive solution of model (\ref{wk12}).
Then the sufficient condition for the distinction of the disease is obtained.
Furthermore, we have obtained the sufficient condition for persistence in mean by  selecting suitable Lyapunov function with regime switching.

Moreover, some interesting topics are deserved for further consideration.
It is  found that  $\widetilde{R}_s^*\leq  {R}_s^*$ in this paper,  and we  will continue to investigate what happens to the diseases under the condition $\widetilde{R}_s^*<1<  {R}_s^*$.
In addition, we can use the similar methods to study more complex epidemic models, such as SEQIR model  with media effects.

\section*{Conflict of interest}
The authors declare that they have no conflict of interest.

\section*{Data availability}
No data is used for the research in this article.

\newpage
\newpage
\bf{Figures}
\begin{figure}[htpb]
  \centering
  {\includegraphics[width=0.6\textwidth,height=0.4\textwidth]{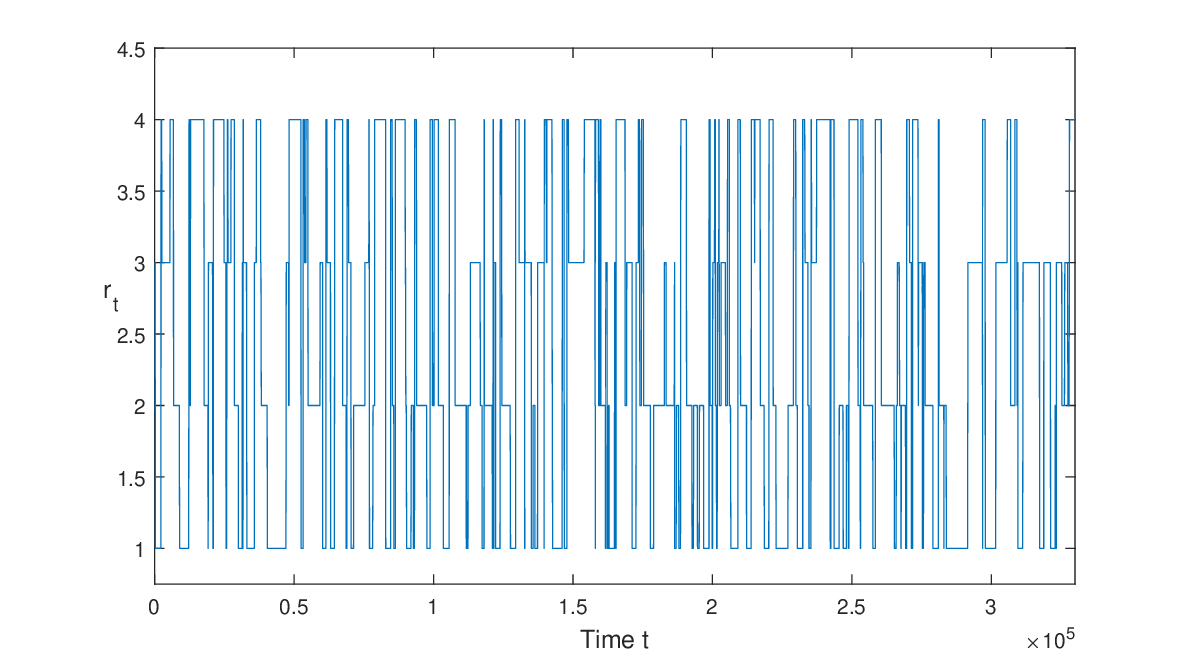}}
\begin{center}
\caption { Simulation of Markov chain $r(t)$ with $r(0)=3$.}
\end{center}
\end{figure}
\begin{figure}[htpb]
  \centering
  {\includegraphics[width=0.7\textwidth,height=0.4\textwidth]{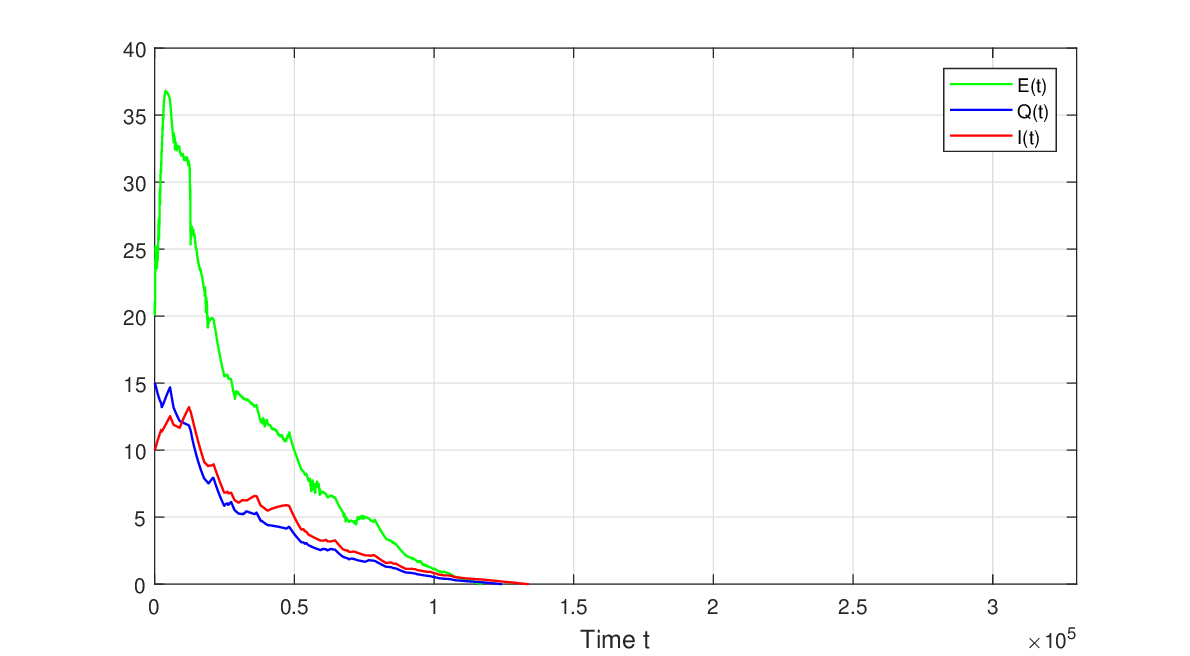}}
\begin{center}
\caption {Simulations of $E(t),Q(t),I(t)$ with $ {R}_0^*<1$.}
\end{center}
\end{figure}
\begin{figure}[htpb]
  \centering
  {\includegraphics[width=0.9\textwidth,height=0.5\textwidth]{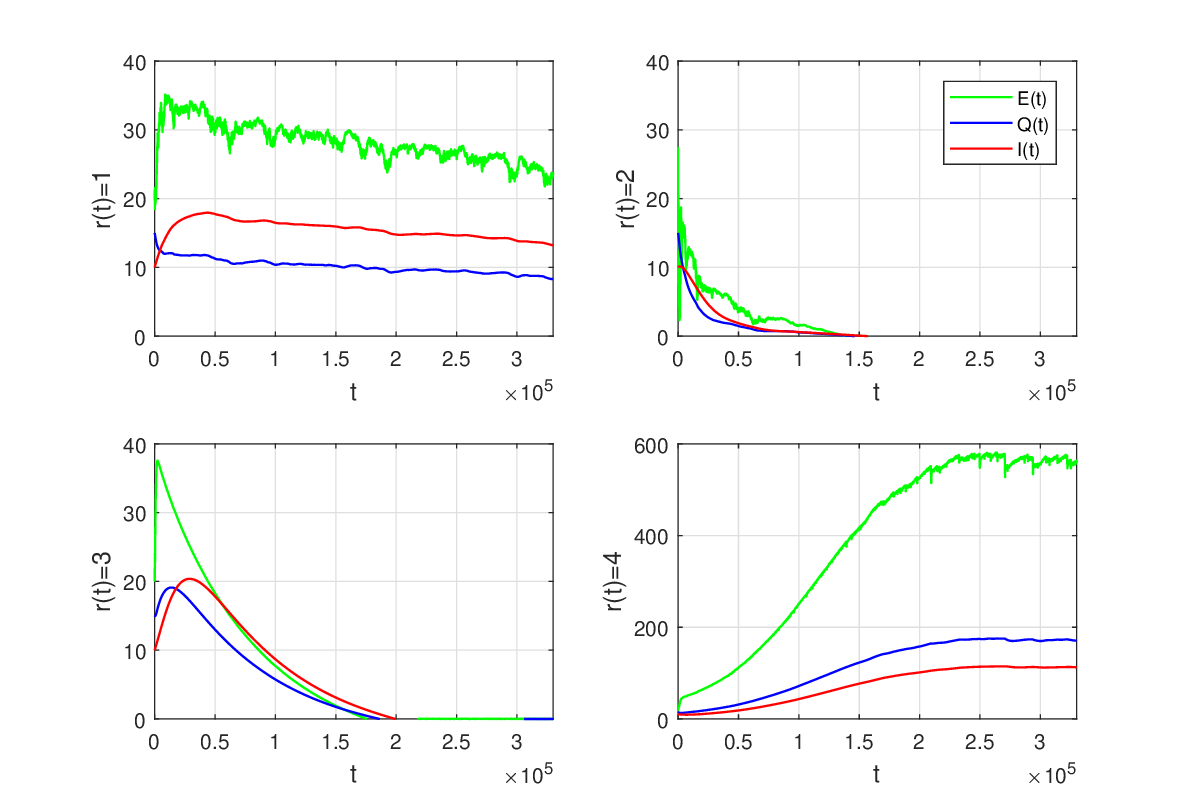}}
\begin{center}
\caption { Simulations of $E(t),Q(t),I(t)$ with $ {R}_0^*<1$  and  $ r(t)=1,2,3,4$.}
\end{center}
\end{figure}
\begin{figure}[htpb]
  \centering
  {\includegraphics[width=0.7\textwidth,height=0.4\textwidth]{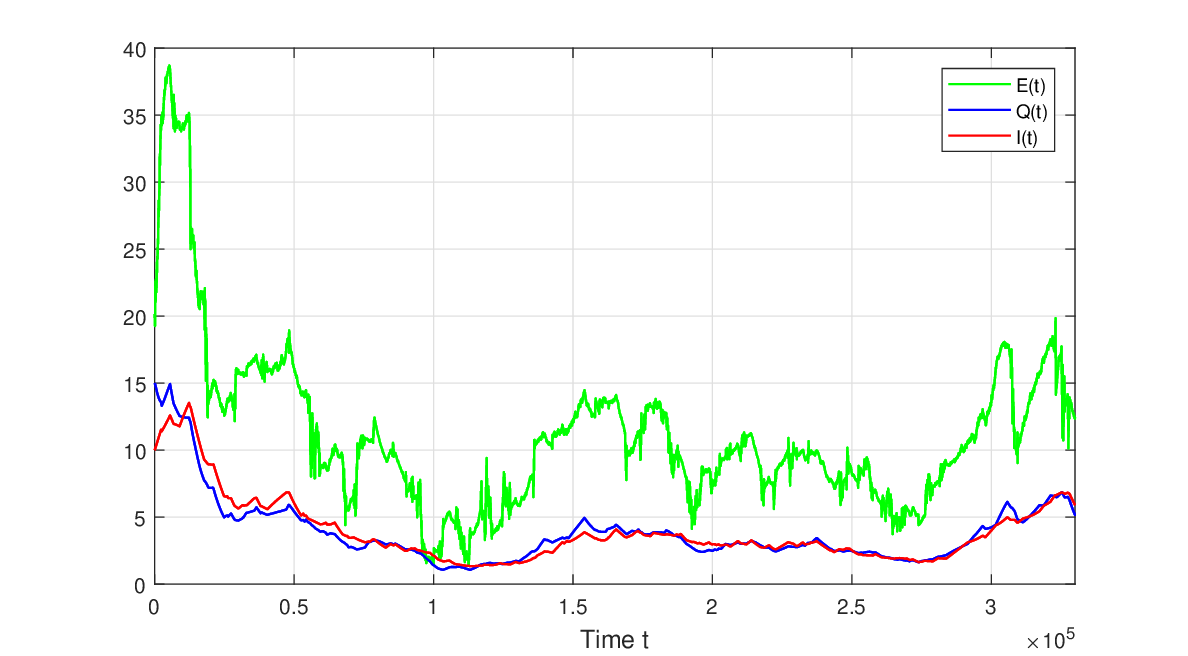}}
\begin{center}
\caption { Simulations of $E(t),Q(t),I(t)$ with   $\widetilde{R}_0^*>1$.}
\end{center}
\end{figure}
\begin{figure}[htpb]
  \centering
  {\includegraphics[width=0.9\textwidth,height=0.5\textwidth]{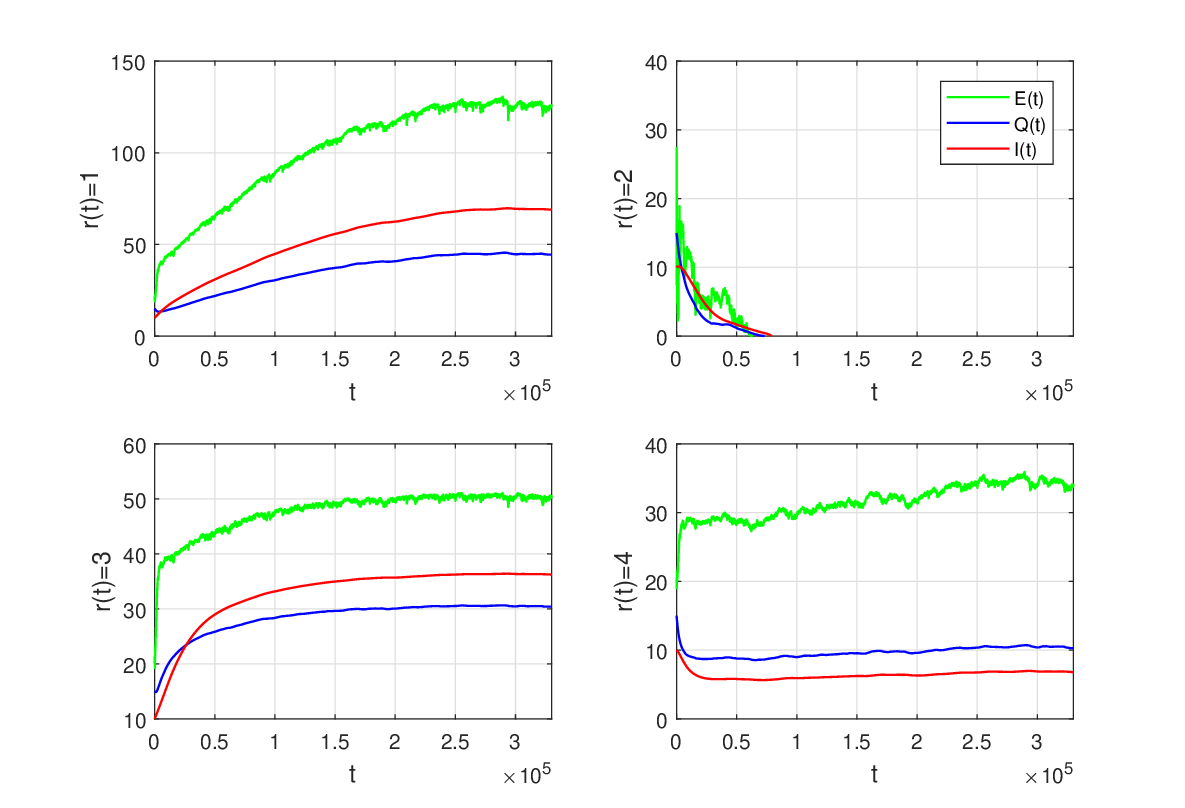}}
\begin{center}
\caption { Simulations of $E(t),Q(t),I(t)$ with $\widetilde{R}_0^*>1$ and  $ r(t)=1,2,3,4$.}
\end{center}
\end{figure}

\end{document}